\documentclass[runningheads,a4paper]{llncs}
\usepackage{makeidx}  
\usepackage{amssymb}
\usepackage{amsmath}
\usepackage{graphicx}
\begin{document}
\mainmatter              
\title{Mimimizing Total Completion Time in Multiprocessor Job Systems with Energy Constraint}
\titlerunning{Mimimizing Total Completion Time under Energy Constraint}  
%
\author{Alexander Kononov\orcidID{0000-0001-6144-0251} \and
Yulia Kovalenko\orcidID{0000-0003-4791-7011}
}
\authorrunning{A. Kononov and Yu. Kovalenko}

%
%
%
\institute{Sobolev Institute of Mathematics, Novosibirsk, Russia\\
\email{alvenko@math.nsc.ru, julia.kovalenko.ya@yandex.ru}
}

\maketitle              

\begin{abstract}
We consider the problem of scheduling multiprocessor jobs to minimize the total completion time under the given energy budget.
Each multiprocessor job requires more than one processor at the same moment of time.
Processors may operate at variable speeds.
Running a job at a slower speed is more energy efficient, however it takes longer time and affects the performance.
The complexity of both parallel and dedicated versions of the problem is investigated.
We propose approximation algorithms for various particular cases.
In our algorithms, initially a sequence of jobs and their processing times are calculated and then a feasible solution is constructed using list-type scheduling rule.
\end{abstract}

\\

\section{Introduction}
We investigate the problem of non-preemptive scheduling a set of jobs $\mathcal{J}=\{1,\dots,n\}$ on $m$ speed scalable parallel processors.
Each job~$j\in \mathcal{J}$ is characterized by processing volume (work) $V_j$ and the number $size_j$ or the set $fix_j$ of required processors.
Note that parameter $size_j$ for job $j \in \mathcal{J}$ indicates that the job can be processed
on any subset of parallel processors of the given size. Such jobs are called rigid jobs~\cite{DrozdowskiBook}.
Parameter $fix_j$ states that the job uses the prespecified subset of dedicated processors. Such jobs are called single mode multiprocessor jobs~\cite{DrozdowskiBook}.
We also consider moldable jobs~\cite{DrozdowskiBook}.
In contrast to the previous job types, a moldable job $j$ may be performed on any number of processors lower or equal to the given upper bound $\delta_j$.

The standard homogeneous model in speed-scaling is considered. When a processor
runs at a speed $s$, then the rate with which the energy is consumed (the \textit{power}) is $s^\alpha$,
where $\alpha > 1$ is a constant (usually, $\alpha \approx 3$).
Each of $m$ processors may operate at variable speed.
However, we assume that the total work $V_j$ of a job $j\in \mathcal{J}$ should be uniformly divided between the utilized processors, i.e.
if job $j$ uses $m_j$ processors, then processing volumes are the same for all $m_j$ processors (denoted by $W_j:=\frac{V_j}{m_j}$), and all these processors run at the same speed.
It is supposed that a continuous spectrum of processor speeds is available.

The aim is to find a feasible schedule with the minimum sum of completion times $\sum C_j$ so that the energy consumption
is not greater than a given energy budget~$E$.
This is a natural assumption in the case when the energy of a battery is fixed, i.e.
the problem finds applications in computer devices whose lifetime depends on a limited battery efficiency
(for example, multi-core laptops). Moreover,  the bicriteria problems
of minimizing energy consumption and a scheduling metric arise in real practice.
The most obvious approach is to bound one of the objective functions and optimize the other.
The energy of the battery may reasonably be estimated, so
we bound the energy used, and optimize the regular timing criterion.

The non-preemptive rigid, moldable and single-mode variants of the speed-scaling scheduling
subject to bound on energy consumption are denoted by ${P|size_j,energy|\sum C_j}$, ${P|var,\delta_j,energy|\sum C_j}$
and ${P|fix_j,energy|\sum C_j}$, respectively.

\section{Previous Research}

Pruhs~\cite{Pruhs2004} et al. investigated the single-processor problem of minimizing the average flow time of jobs, given a fixed amount of energy and release times of jobs. For unit-work jobs,
they proposed a polynomial-time algorithm that simultaneously computes, for each
possible energy level, the schedule with smallest average flow time. Bunde~\cite{Bunde2009} adopted the approach to multiple processors.
$O(1)$-approximation algorithm, allowing an additional factor of $(1 + \varepsilon)$
energy, has been proposed for scheduling arbitrary work jobs on single processor.
Albers and Fujiwara~\cite{Albers2007} have investigated online and offline versions of single-processor scheduling to minimize energy consumption plus job flow times.
A deterministic constant competitive online algorithm and offline dynamic programming algorithm with polynomial time complexity were proposed for unit-work jobs.

Shabtay et al.~\cite{SK2006} analyzed a closely related problem of scheduling single-processor jobs
on identical parallel processors,
where job-processing times $p_j$ are controllable through the allocation of a nonrenewable common limited resource as
$p_j(R_j)=\left(\frac{W_j}{R_j}\right)^{\kappa}$. Here
$W_j$ is the workload of job $j$, $R_j$ is the amount of resource allocated to processing job $j$
and $0<\kappa \le 1$ is a  positive constant. Exact polynomial time algorithm was proposed
for the multiprocessor non-preemptive instances of minimizing the sum of completion times.
The algorithm can be adopted to the speed scaling scheduling of single-processor jobs.

The speed scaling scheduling with  makespan criterion has been widely investigated.
Various approaches to construct approximation algorithms for single-processor and multiprocessor jobs were proposed (see, e.g., \cite{KonKov2020,Bunde2009,Pruhs2004}).

Now we review the known results for the the classic problem of scheduling multiprocessor jobs
with given durations and without energy constraint. The non-preemptive problem with rigid jobs is strongly NP-hard even in the case of two processors~\cite{Lee1999},
the preemptive one is NP-hard when the number of processors is a part of the input~\cite{Drozdowski2000}.
For the non-preemptive instances constant factor approximation algorithms were proposed. These algorithms use list-type scheduling~\cite{TLWFTGSY1994} and
scheduling to minimize average response time (SMART)~\cite{Schwiegelshohn1998}.
The non-preemptive single-mode problem is NP-hard in the case of two processors~\cite{Hoogeveen1994} and strongly NP-hard in the case of two-processor jobs~\cite{Kubale1996}.
The strategy from preemptive schedule to the non-preemptive one gives a $2$-approximation algorithm for two-processor problem~\cite{Cai1998}, and
First Fit Coloring strategy allows to obtain a $2$-approximate solution for unit-work two-processor jobs~\cite{Giaro1999}.

\subsubsection*{Our results.}
We prove NP-hardness of problems $P|size_j\le \frac{m}{2},energy|\sum C_j$ and $P|fix_j, \ |fix_j|=2,energy|\sum C_j$, and develop
two-stage approximation algorithms for the following particular cases: 
\begin{itemize}
\item rigid jobs requiring at most $\frac{m}{2}$ processors,
\item moldable jobs,
\item two-processor dedicated instances.
\end{itemize}
At the first stage,
we obtain a lower bound on the total completion time and calculate a sequence and processing  times  of  jobs using an auxiliary convex program.
Then, at the second stage, we transform  our  problem  to  the  classic scheduling problem  without  speed  scaling,  and  we
use ''list-scheduling`` algorithms to obtain feasible solutions. Whenever a subset of processors falls idle, a ''list-scheduling`` algorithm schedules
from a given priority list the first job that does not require more processors than are available.

\section{Rigid Jobs}

In this section we consider rigid jobs.
Firstly, we prove that the problem is NP-hard.
Secondly, $2$-approximation algorithm is presented for jobs, which require at most $\frac{m}{2}$ processors
and have identical workloads on utilized processors (i.e. $V_j=W\cdot size_j$).

\subsection{NP-hardness}

\begin{theorem}
\label{th:rigidNP}
Problem $P|size_j\le \frac{m}{2},W_j=1,energy|\sum C_j$ is NP-hard in the strong sense.
\end{theorem}

\textbf{Proof.} We show that the strongly NP-complete 3-PARTITION problem polynomially transforms to
the decision version of scheduling problem $P|size_j\le \frac{m}{2},W_j=1,energy|\sum C_j$.

We consider an instance of the 3-PARTITION problem: Given a set of $3q$ elements with weights $a_j,\ j=1,\dots,3q,$ where
$\sum_{j=1}^{3q} a_j=Bq$ and $\frac{B}{4}\le a_j\le \frac{B}{2}$.
Could the set be partitioned into $q$ subsets $A_1,\dots, A_q$ such that $\sum_{j\in A_i} a_j=B$?

An instance of $P|size_j,energy|\sum C_j$ is constructed as follows.
Put the number of jobs $n=3q$, the number of processors $m=B$, and the energy budget $E=Bq$.
For every $a_j$ we generate a job $j,\ j=1,\dots,3q$.
We set $W_j=1$, $size_j=a_j$, $V_j=a_j$ for $j\in \mathcal{J}$.
In the decision version of $P|size_j,W_j=1,energy|\sum C_j$ it is required to answer the question: Is there
a schedule with $\sum C_j$ value not greater than a given threshold $T$?

In order to determine the value of $T$ we solve an auxiliary problem with $\sum_{j=1}^{3q} a_j$  single-processor jobs of unit works,
i.e. each rigid job is replaced by $size_j$ single-processor jobs.
Such problem has the unique optimal solution (with the accuracy of placing jobs on processors and permuting jobs on each processor),
where each processor executes $q$ jobs and uses energy budget $q$.
Now we find optimal durations of jobs on each processor, solving the following convex program:

\begin{equation}
\sum_{j=1}^q p_j(q-j+1)\to \min,
\end{equation}

\begin{equation}
\sum_{j=1}^q p_j^{1-\alpha}=q,
\end{equation}

\begin{equation}
p_j\ge 0,\ j=1,\dots,q.
\end{equation}

Here $p_j$ is the execution time of $j$-th job on a processor, $j=1,\dots,q$.

We compose the Lagrangian function $L(p_j,\lambda)=\sum_{j=1}^q p_j(q-j+1)+$\linebreak
$\lambda\left( \sum_{j=1}^q p_j^{1-\alpha}-q \right)$ and calculate the optimal solution
by equating to zero partial derivatives:

$$
p_j^{*}=\frac{\left(\sum_{j=1}^q (n-j+1)^{\frac{\alpha-1}{\alpha}}\right)^{\frac{1}{\alpha-1}}} {q^{\frac{1}{\alpha-1}} (q-j+1)^{\frac{1}{\alpha}}},\ j=1,\dots,q,
$$

$$
\sum C_j^{*}=\sum_{j=1}^q p_j(q-j+1)=\left(\sum_{j=1}^q (n-j+1)^{\frac{\alpha-1}{\alpha}}\right)^{\frac{\alpha}{\alpha-1}} q^{\frac{1}{1-\alpha}}.
$$

Note that each next job has more duration than the previous one.
The optimal schedule for each processor does not have idle times.
The optimal total completion time for all processors is equal to $m \sum C_j^{*}$.

Set the threshold $T:=3 \sum C_j^{*}$, since at most 3 rigid jobs can be executed in parallel. 
We show that a positive answer (a negative answer) to 3-PARTITION implies a positive answer (a negative answer)
to the constructed $P|size_j,energy|\sum C_j$ with $\sum C_j\le T$.

Firstly, we assume that the answer to 3-Partition is positive. Then there is a feasible schedule with  $\sum C_j\le T$, where
three jobs, corresponding to three elements forming set $A_i,\ i=1,\dots, q$ such that $\sum_{j\in A_i} a_j=B$,
are executed in parallel. This schedule is similar to the optimal schedule of the corresponding problem with single-processor jobs.
The value of criterion is equal to $3 \sum C_j^{*}$. 

Secondly, we show that the negative answer to 3-Partition implies the negative answer to our speed scaling scheduling problem.
Indeed, in this case we can not construct a schedule, which is identical to the optimal schedule for the corresponding single-processor jobs
(with the accuracy of placing jobs on processors and permuting jobs on each processor).
In other words, in any feasible schedule the sum of completion times is greater than $3 \sum C_j^{*}$  as this schedule has idle times.

The presented transformation is polynomial. So, problem $P|size_j\le \frac{m}{2},W_j=1,energy|\sum C_j$ is strongly NP-hard. \qed

In the next subsection, we present a polynomial time approximation algorithm with constant factor approximation guarantee.
Processing times of jobs are calculated using convex program and approximate schedule is constructed by ``list-scheduling'' algorithm.

\subsection{Approximation Algorithm} 
\label{subsec:LBRigid}
The sequence of jobs and the completion time of each job
are important in the problems  with criterion~$\sum C_j$.
Now we compute a lower bound for the case when a jobs sequence is given.
Suppose that the jobs are started in accordance with permutation $\pi=(\pi_1,\dots,\pi_n)$. 
Using the lower bound on the total completion time presented in~\cite{TLWFTGSY1994} for rigid jobs with the given durations, we formulate the following convex program:

\begin{equation}
\label{LB_Csum}
\frac{1}{m}\sum_{j=1}^n\sum_{i=1}^{j}size_{\pi_i}p_{\pi_i} +\frac{1}{2}\sum_{j=1}^np_{\pi_j}-
\frac{1}{2m}\sum_{j=1}^{n}size_{\pi_j}p_{\pi_j}\to \min,
\end{equation}
\begin{equation}
\label{LB_Energy}
\sum_{j=1}^n W_j^{\alpha} p_j^{1-\alpha} size_j \le E.
\end{equation}
\begin{equation}
\label{LB_Def}
p_j\ge 0,\ j\in \mathcal{J}.
\end{equation}

We solve the program by means of the Lagrangian method.
Define the Lagrangian function $L(p_{\pi_j},\lambda)$ as

$$L(p_{\pi_j},\lambda)=\frac{1}{m}\sum_{j=1}^n\sum_{i=1}^{j}  size_{\pi_i}p_{\pi_i} +\frac{1}{2}\sum_{j=1}^n\left( 1- \frac{size_{\pi_j}}{m} \right)p_{\pi_j}+$$
$$\lambda \left( \sum_{j=1}^n W_{\pi_j}^{\alpha} p_{\pi_j}^{1-\alpha} size_{\pi_j} - E \right).$$

The necessary and sufficient conditions for an optimal solution are (partial derivatives are equal to zero):

$$\frac{\partial L}{\partial p_{\pi_i}}=
\frac{1}{m} size_{\pi_i} (n-i+1) +\frac{1}{2}\left(1-\frac{size_{\pi_i}}{m}\right)
+ \lambda W_{\pi_i}^{\alpha} size_{\pi_i} (1-\alpha) p_{\pi_i}^{-\alpha} =0,$$
$$i=1,\dots,n.$$

Rewriting the expressions, we obtain

$$p_{\pi_i}=\left((\alpha-1)\lambda m\right)^{\frac{1}{\alpha}}
\frac{W_{\pi_i} size_{\pi_i}^{\frac{1}{\alpha}}}{\left(size_{\pi_i}(n-i+0.5)+0.5m\right)^{\frac{1}{\alpha}}},$$
$$i=1,\dots,n.$$

The processing times are placed into equation

$$\frac{\partial L}{\partial \lambda}=
 \sum_{j=1}^n W_j^{\alpha} p_j^{1-\alpha} size_j - E=0.$$

As a result we calculate the durations of jobs

$$
p_{\pi_i}=\frac{E^{\frac{1}{1-\alpha}} W_{\pi_i} size_{\pi_i}^{\frac{1}{\alpha}}}
{\left(size_{\pi_i}(n-i+0.5)+0.5m\right)^{\frac{1}{\alpha}}}
\cdot
\left(
\sum_{j=1}^n
\frac{W_{\pi_j} size_{\pi_j}^{\frac{1}{\alpha}}}
{\left(size_{\pi_j}(n-j+0.5)+0.5m\right)^{\frac{1-\alpha}{\alpha}}}
\right)^{\frac{1}{\alpha-1}},
$$
$$i=1,\dots,n.$$

The obtained values for execution times are placed in expression (\ref{LB_Csum})
and the lower bound on $\sum C_j$ in the general case for an arbitrary permutation $\pi$
of jobs is calculated  as follows
\begin{equation}
\label{F:permut}
LB(\pi)=\frac{E^{\frac{1}{1-\alpha}}}{m}
\left(
\sum_{i=1}^n W_{\pi_i} size_{\pi_i}^{\frac{1}{\alpha}}
\left( size_{\pi_i}(n-i+0.5)+0.5m\right)^{\frac{\alpha-1}{\alpha}} \right)^{\frac{\alpha}{\alpha-1}}.
\end{equation}
Here $\pi_i$  is the $i$-th job in accordance with permutation $\pi$.
So, it is required to find permutation, that gives $\min_{\pi}LB(\pi)$.

From now on, we suppose that the processing works of jobs on processors are identical, i.e. $W_j=W,\ j\in \mathcal{J}$.
Then the minimization of (\ref{F:permut}) is equivalent to the minimization of

$$
G(\pi)=
\sum_{i=1}^n W size_{\pi_i}
\left( n-i+0.5+\frac{0.5m}{size_{\pi_i}}\right)^{\frac{\alpha-1}{\alpha}}.
$$

We define vectors $WS_{\pi}=(W size_{\pi_1},W size_{\pi_2},\dots,W size_{\pi_n})$ and
$NS_{\pi}=\left(\left(n-0.5+\frac{0.5m}{size_{\pi_1}}\right)^{\frac{\alpha-1}{\alpha}},\left(n-1.5+\frac{0.5m}{size_{\pi_2}}\right)^{\frac{\alpha-1}{\alpha}},\dots,
\left(0.5+\frac{0.5m}{size_{\pi_n}}\right)^{\frac{\alpha-1}{\alpha}}\right)$.
Then $G(\pi)$ can be expressed as the following scalar product $WS_{\pi} \bullet (NS_{\pi})^{\mathrm{T}}$.
It is easy to see that the minimum of $WS_{\pi} \bullet (NS_{\pi})^{\mathrm{T}}$ is reached on the permutation,
where the jobs are ordered by non-decreasing of the required processors numbers $size_j,\ j\in \mathcal{J}$.
Indeed, in this case the first vector is ordered nondecreasingly, the second one is ordered nonincreasingly
and the scalar product is minimal (the proof is based on the permutation approach).

Let the jobs are ordered by non-decreasing of $size_j$. We denote by $\bar{p}_j$  the durations of jobs, and by
$LB(\bar{p}_j)=\frac{1}{m}\sum_{j=1}^n\sum_{i=1}^{j-1}size_j\bar{p}_j +\frac{1}{2}\sum_{j=1}^n\bar{p}_j+$\linebreak
$\frac{1}{2m}\sum_{j=1}^{n}size_j\bar{p}_j$ the lower bound corresponding to the optimal solution of problem~(\ref{LB_Csum})-(\ref{LB_Def}).
We construct the schedule using 
``list-scheduling'' algorithm with processing times $\bar{p}_j,\ j\in \mathcal{J}$: The first job is scheduled at time $0$.
The next job is scheduled at the earliest time such that there are enough processors to execute it.
Recall that each job requires at most $\frac{m}{2}$ processors and $W_j=W$. 

``List-scheduling'' algorithm has time complexity $O(n^2)$ and
allows to construct $2$-approximate schedule.
Indeed, the starting time of the $j$-th job is not greater than $2\left(\frac{1}{m} \sum_{i=1}^{j-1}size_i\bar{p}_i \right)$,
as at least $\frac{m}{2}$ processors are busy at each time moment in the schedule.
The completion time of the $j$-th job satisfies condition $C_j\le \frac{2}{m} \sum_{i=1}^{j-1}size_i\bar{p}_i+\bar{p}_j$.
In the sum we have
$$\sum_{j\in \mathcal{J}}C_j\le
2\left(\frac{1}{m}\sum_{j=1}^{n}\sum_{i=1}^{j-1} size_i\bar{p}_i+\frac{1}{2}\sum_{j=1}^{n}\bar{p}_j\right).$$
We compare this value with the lower bound $LB(\bar{p}_j)$ 
and conclude that $\sum_{j\in \mathcal{J}}C_j\le 2 LB(\bar{p}_j)$.

Therefore, the following theorem takes place.

\begin{theorem}
A $2$-approximate schedule
can be found in polynomial time for problem $P|size_j,W_j=W,size_j\le \frac{m}{2},energy|\sum C_j$.
\end{theorem}

\section{Moldable Jobs}
\label{subsec:CSMoldable}
Now we provide $2$-approximation algorithm for moldable jobs with identical works.
Recall that $V_j$ denote the total processing volume of job $j$, i.e. the execution time of this job on one processor with unit speed.
Let $\delta_j\le m$ be the maximal possible number of processors, that may be utilized by job~$j$.

In order to obtain a lower bound on the sum of completion times,
we formulate the following convex model in the case of the given sequence $\pi$ of jobs
\begin{equation}
\label{LB_CsumMoldable}
\frac{1}{m}\sum_{j=1}^n\sum_{i=1}^{j}p_{\pi_i} +\frac{1}{2}\sum_{j=1}^n\frac{p_{\pi_j}}{\delta_{\pi_j}}-
\frac{1}{2m}\sum_{j=1}^{n}p_{\pi_j}\to \min,
\end{equation}

\begin{equation}
\sum_{j=1}^n V_j^{\alpha} p_j^{1-\alpha}\le E.
\end{equation}
Here $p_j$ is the execution time of job~$j$ on one processor in the total volume $V_j$.

Using the same arguments as in Subsection~\ref{subsec:LBRigid},
we can show that the minimum of (\ref{LB_CsumMoldable}), say
$$LB:= \frac{E^{1/1-\alpha}}{m} \left(\sum_{j=1}^{n} V_{\pi_j} \left( n-j+0.5 +\frac{0.5m}{\delta_{\pi_j}} \right)^{\alpha-1/\alpha}  \right)^{\alpha/\alpha-1},$$
is reached on the permutation $\pi$,
where the jobs are ordered by non-decreasing of the maximum processors numbers $\delta_j,\ j\in \mathcal{J}$,
in the case when works of jobs $V_j$ are identical or the non-decreasing order of $V_j$
corresponds to the non-decreasing order of $\delta_j$ (i.e. $V_i<V_j$ implies $\delta_i\le \delta_j$).
The corresponding durations of jobs will be denoted by $\bar{p}_j,\ j\in \mathcal{J}$.

Now we assign the number of processors $m_j$ for jobs as follows
\[ m_j = \begin{cases}  \delta_j     & \text{if $\delta_j<\left\lceil \frac{m}{2} \right\rceil$,} \\
   \left\lceil \frac{m}{2} \right\rceil  & \text{if $\delta_j\ge \left\lceil \frac{m}{2} \right\rceil$,}
\end{cases} \]
and construct the schedule using ``list-scheduling'' algorithm based on the order of jobs
in non-decreasing of $\delta_j,\ j\in \mathcal{J}$. Let us prove that the total completion time $\sum C_j(\bar{p}_j)\le 2 LB$.
Indeed,

$$
\sum C_j(\bar{p}_j)\le \frac{2}{m} \sum_{j=1}^n\sum_{i=1}^{j-1} \bar{p}_i +\sum_{j=1}^n \frac{\bar{p}_j}{m_j}=
\frac{2}{m} \sum_{j=1}^n\sum_{i=1}^{j} \bar{p}_i +\sum_{j=1}^n \frac{\bar{p}_j}{m_j}-\frac{2}{m}\sum_{j=1}^n \bar{p}_j$$
$$=
\frac{2}{m} \sum_{j=1}^n\sum_{i=1}^{j} \bar{p}_i +  \sum_{j=1}^n \bar{p}_j \left(\frac{1}{m_j}-\frac{2}{m}\right)\le
\frac{2}{m} \sum_{j=1}^n\sum_{i=1}^{j} \bar{p}_i + \sum_{j=1}^n \bar{p}_j \left(\frac{1}{\delta_j}-\frac{1}{m}\right)\le
2LB.
$$

Therefore, we have

\begin{theorem}
A $2$-approximate schedule
can be found in polynomial time for problem $P|any,V_j=V,\delta_j,energy|\sum C_j$.
\end{theorem}

We note here, that the complexity status of speed scaling scheduling problem $P|any,V_j=V,\delta_j,energy|\sum C_j$ is open.

\section{Single Mode Multiprocessor Jobs}

In this section we consider single-mode multiprocessor jobs.
Firstly, we prove that the problem is NP-hard. Then we provide a polynomial time algorithm for two-processor instances.

\subsection{NP-hardness}

\begin{theorem}
Problem $P|fix_j,\ |fix_j|=2,V_j=2,energy|\sum C_j$ is NP-hard in the strong sense.
\end{theorem}

\textbf{Proof.}
The proof is similar to the proof of Theorem~\ref{th:rigidNP}, but it is based on the polynomial reduction
of the strongly NP-complete Chromatic Index problem for cubic graphs~\cite{Holyer1981}.

The chromatic index of a graph is the minimum number of colors
required to color the edges of the graph in such a way that no two adjacent
edges have the same color. Consider an instance of the Chromatic Index problem on a cubic graph $G=(V, A)$,
which asks whether the chromatic index $\chi'(G)$ is three.
It is well-know that $\chi'(G)= 3$ or $4$, and $|V|$ is even.
Moreover, $\chi'(G)= 3$ if and only if each color class has exactly $\frac{1}{2}|V|$ edges.

We construct an instance of $P|fix_j,energy|\sum C_j$ as follows.
Put the number of jobs $n=|A|$, the number of processors $m=|V|$ and the energy budget $E=2|A|=3|V|$.
Vertices correspond to processors.
For every edge $\{u_j,v_j\}$ we generate a job $j$ with $fix_j=\{u_j,v_j\}$, $j=1,\dots,|A|$.
We set $V_j=2$ and $W_j=\frac{V_j}{|fix_j|}=1$ for all $j=1,\dots,n$.
In the decision version of  $P|fix_j,W_j=1,energy|\sum C_j$ it is required to answer the question: Is there
a schedule, in which the total completion time is not greater than a given threshold~$T$?

In order to define the value of $T$ we solve auxiliary problem with $2n$  single-processor jobs,
i.e. each two processor job is replaced by two single-processor jobs.
Such problem has the unique optimal solution (with the accuracy of permuting jobs on processors),
where each processor execute tree jobs and uses energy budget $3$.
Now we find optimal durations of jobs on each processor, solving the following convex program:

\begin{equation} 
p_1 + 2p_2 +p_3\to \min,
\end{equation}

\begin{equation}
\sum_{j=1}^3 p_j^{1-\alpha}=3,
\end{equation}

\begin{equation}
p_j\ge 0,\ j=1,\dots,n.
\end{equation}

Here $p_j$ is the execution time of $j$-th job on a processor.

We compose the Lagrangian function
$$L(p_j,\lambda)=(3p_1 + 2p_2 +p_3)+ \lambda\left( p_1^{1-\alpha}+p_2^{1-\alpha}+p_3^{1-\alpha}-3 \right)$$
and calculate

\begin{equation}
p^{*}_j=\frac{3^{1/1-\alpha}}{(4-j)^{1/\alpha}} \left( 3^{\alpha-1/\alpha} + 2^{\alpha-1/\alpha} + 1\right)^{1/\alpha-1},\ j=1,2,3,
\end{equation}
\begin{equation}
\sum C_j^{*}=\left( \frac{\left( 3^{\alpha-1/\alpha} + 2^{\alpha-1/\alpha} + 1\right)^{\alpha} }{3} \right)^{1/\alpha-1}.
\end{equation}

The sum of job completion times on all processors is equal to $m \sum C_j^{*}$. The optimal schedule does not have idle times.
Set the threshold $T:=\frac{m}{2} \sum C_j^{*}$ because at most $\frac{m}{2}$ two processor jobs can be executed in parallel.
We prove that a positive answer (a negative answer) to Chromatic Index problem corresponds a positive answer (a negative answer)
to the constructed decision version of $P|fix_j,W_j=1,energy|\sum C_j$.

Now we assume that the answer to Chromatic Index problem is positive. Then there is a feasible schedule, where 
$\frac{m}{2}$ jobs, corresponding to $\frac{m}{2}$ edges forming one coloring class,
are executed in parallel. This schedule is similar to the optimal schedule of the corresponding problem with single-processor jobs.
The value of criterion is equal to $\frac{m}{2} \sum C_j^{*}$. 

It is easy to see that the negative answer to Chromatic Index problem implies the negative answer to our scheduling problem.
Indeed, in this case we can not construct a schedule, which is identical to the optimal schedule for the corresponding single-processor jobs
(with the accuracy of permuting jobs on processors).
In other words, any feasible schedule has idle times, and, therefore, the total sum of completion times is greater than $\frac{m}{2} \sum C_j^{*}$. 

The presented reduction is polynomial. So, speed scaling scheduling problem $P|fix_j,W_j=1,energy|\sum C_j$ is strongly NP-hard. \qed

\subsection{Two-processor Instances}
Now we consider the non-preemptive problem with two processors and propose polynomial time algorithm with constant-factor approximation guarantee.
Denote by $\mathcal{J}_i$ the set of jobs using only processor $i=1,2$
and by $\mathcal{J}_{12}$ the set of two-processors jobs, $\mathcal{J}=\mathcal{J}_1\cup \mathcal{J}_2\cup \mathcal{J}_{12}$.
In order to obtain a lower bound on the total completion time we identify two subproblems: the first one schedules only jobs from $\mathcal{J}'=\mathcal{J}_1\cup \mathcal{J}_{12},\ |J'|=n'$,
the second one schedules only jobs from $\mathcal{J}''=\mathcal{J}_2\cup \mathcal{J}_{12},\ |J''|=n''$.

The optimal solution of the first subproblem in the case of the given sequence $\pi$ of jobs can be found by solving the following convex program:
$$\sum_{i=1}^{n'} (n'-i+1)p_{\pi_i}\to \min,$$

$$\sum_{i \in \mathcal{J}'} |fix_i| {(p_i)}^{1-\alpha}W_i^{\alpha}\le E.$$

Solving this subproblem via KKT-conditions, we obtain durations of jobs
$$p_{\pi_i}=\left( \frac{E}{\sum_{j=1}^{n'} W_{\pi_j} |fix_{\pi_j}|^{1/\alpha} (n'-j+1)^{\alpha-1/\alpha}}\right)^{1/1-\alpha} \frac{W_{\pi_i} |fix_{\pi_i}|^{1/\alpha} }{(n'-i+1)^{1/\alpha}},
$$
$$i=1,\dots,n',$$
and the sum of completion times
$$\sum C_j^1 (\pi)=\left(E\right)^{1/1-\alpha} \left( \sum_{j=1}^{n'} W_{\pi_j} |fix_{\pi_j}|^{1/\alpha} (n'-j+1)^{\alpha-1/\alpha}\right)^{\alpha/\alpha-1}.$$
So, the minimum sum of completion times is reached on the permutation $\pi'$,
where the jobs are ordered by non-decreasing of $W_i|fix_i|^{1/\alpha},\ i\in \mathcal{J}'$,
since values $(n'-j+1)$ decrease.

Using the same approach for the second subproblem we conclude that the minimum sum of completion times $\sum C_j^2$ is reached on the permutation $\pi''$,
where the jobs are ordered by non-decreasing of $W_i|fix_i|^{1/\alpha},\ i\in \mathcal{J}''$.
So, the lower bound $LB$ for the general problem can be calculated as\linebreak $\max\{\sum C_j^1 (\pi'),\ \sum C_j^2 (\pi'')\}$.
 Note that subsequences of two-processor jobs are identical in optimal solutions of both subproblems.

Decreasing the energy budget in both subproblems in two times,
we obtain $2^{1/\alpha-1}$-approximate solutions $S'$ and $S''$ for them with the same sequences of jobs as in case of energy budget $E$.
Let $C_j'$ and $C_j''$ ($p_j'$ and $p_j''$) denote the completion times (the processing times) of two processor job $j$ in $S'$ and $S''$, respectively.
Now we construct a preemptive schedule for the general problem, and then transform this schedule to the non-preemptive one.

In the constructed preemptive schedule each two-processor job $j$ is executed without preemptions
in interval $( \max\{C_j', C_j''\}-\min\{p_j', p_j''\},\ \max\{C_j', C_j''\}]$.
Note that execution intervals of two-processor jobs do not intersect each other.
Single-processor jobs are performed in the same order and with the same durations as in $S'$ and $S''$,
but may be preempted by two-processor jobs (idle times between single-processor jobs are not allowed).
It is easy to see that the schedule is feasible,
and has the total completion time $\sum C_j\le 2\cdot 2^{1/\alpha-1} LB$ as the completion times of single-processor jobs are no later than in $S'$ and $S''$ by construction.

Now we go to calculate a non-preemptive feasible schedule.
The obtained preemptive schedule may be reconstructed without increasing the completion times of jobs such that at most one single-processor job is preempted by
each two-processor job (if a two-processor job $j$ preempts two single-processor jobs, then moving $j$ slightly earlier will lower the completion time of $j$  without
affecting the completion times of any other jobs). Let $S'(j)$ and $C'(j)$  denote the starting time and completion time, respectively, of any job $j$
in the reconstructed preemptive schedule.
Identify single-processor jobs $j_{i_1},\ j_{i_2},\dots,j_{i_k}$ that are preempted by some two-processor jobs.
Suppose that these jobs are ordered by increasing of starting times $S'(j_{i_1})<\dots<S'(j_{i_k})$, and therefore completion times satisfy $C'(j_{i_1})<\dots<C'(j_{i_k})$.
Let $F(j_{i_l})$  be  the last two-processor job that preempts $j_{i_l}$,
$g(j_{i_l})$ be the amount of processing time of $j_{i_l}$ scheduled before the starting time of $F(j_{i_l})$,
$h(j_{i_l})$ be the number of jobs that complete later than $F(j_{i_l})$, $l=1,\dots,k$.
The construction procedure consists of $k$ steps.
At step $l$ we insert an idle time period of length $g(j_{i_l})$ on both processors immediately after the
completion of $F(j_{i_l})$. Change the start time of $j_{i_l}$ to the completion time of $F(j_{i_l})$.
So, at step $l$, inserting the idle time period will increase the total
completion time of the schedule by $g(j_{i_l})\cdot h(j_{i_l})$, and in total after all steps the non-preemptive schedule has the objective value
$$\sum C_j^{npr}\le \sum C_j + \sum_{l=1}^{k} g(j_{i_l})\cdot h(j_{i_l}).$$
Since each two-processor job preempts at most one single-processor job,
the time intervals $(S'(j_{i_1}),C'(F(j_{i_1})]$, $(S'(j_{i_2}),C'(F(j_{i_2})]$,$\dots$, $(S'(j_{i_k}),C'(F(j_{i_k})]$
do not overlap with each other. Therefore,
$$\sum C_j \ge \sum_{l=1}^{k} \left(C'(F(j_{i_l})) - S'(j_{i_l})\right)\cdot h(j_{i_l})>\sum_{l=1}^{k} g(j_{i_l})\cdot h(j_{i_l}).$$

As a result we have the following bound on the total completion time
$$\sum C_j^{npr}\le \sum C_j + \sum_{l=1}^{k} g(j_{i_l})\cdot h(j_{i_l})<2 \sum C_j\le 2^{2\alpha-1/\alpha-1} LB.$$
\begin{theorem}
A $2^{2\alpha-1/\alpha-1}$-approximate schedule
can be found in polynomial time for problem $P2|fix_j,energy|\sum C_j$.
\end{theorem}
\begin{corollary}
A $2^{\alpha/\alpha-1}$-approximate schedule
can be found in polynomial time for problem $P2|fix_j,pmtn,energy|\sum C_j$.
\end{corollary}
The complexity status of both preemptive and non-preemptive problems with two processors is open.

\section*{Conclusion}
NP-hardness of both parallel and dedicated versions
of the speed scaling problem with the total completion time criterion under the given energy budget is proved.
We propose an approach to construct approximation algorithms for various particular cases of the problem.
In our algorithms, initially a sequence of jobs and their processing times are calculated and then a feasible solution is constructed using list-type scheduling rule.

Further research might address the approaches to the problems with more complex
structure, where processors are heterogeneous and jobs have alternative execution
modes with various characteristics. Open questions are the complexity status of the problem with moldable jobs and two-processor dedicated problem, and
constant factor approximation guarantee for two-processor jobs in the system with arbitrary number of processors.

\subsubsection*{Acknowledgements}
The reported study was funded by RFBR, project number 20-07-00458.

\end{document}